\documentclass[11pt,a4,fleqn]{article}
\usepackage{graphicx}
\usepackage{amsmath,amssymb,latexsym,graphics,epsfig}
\usepackage{hyperref}
\usepackage{color}
\usepackage{amsthm}

\newtheorem{theorem}{\bf Theorem}[section]
\newtheorem{proposition}[theorem]{\bf Proposition}
\newtheorem{lemma}[theorem]{\bf Lemma}
\newtheorem{corollary}[theorem]{\bf Corollary}

\newtheorem*{observation}{\bf Observation}

\theoremstyle{definition}

\newtheorem{example}{Example}

\newcommand{\lal}{\lambda_{\ell}}
\newcommand{\mul}{\mu_{\ell}}
\newcommand{\nul}{\nu_{\ell}}
\newcommand{\G}{\Gamma}

\newcommand{\pf}{\noindent{\em Proof: }}
\newcommand{\epf}{\hfill\hbox{\rule{3pt}{6pt}}\\}

\def\sgn{\mathop{\rm sign }\nolimits}

\numberwithin{equation}{section}

\begin{document}
\title{{\Large Directed strongly walk-regular graphs}}

\author{E.R. van Dam$^{\textrm{a}}$   \quad  G.R. Omidi$^{\textrm{b},\textrm{c},1}$  \\[2pt]
{\small  $^{\textrm{a}}$Tilburg University, Dept. Econometrics and Operations Research,}\\
{\small  P.O. Box 90153, 5000\,LE, Tilburg, The Netherlands}\\
{\small  $^{\textrm{b}}$Dept. Mathematical Sciences, Isfahan University of Technology},\\
{\small Isfahan, 84156-83111, Iran}\\
{\small $^{\textrm{c}}$School of Mathematics, Institute for Research in Fundamental Sciences (IPM),}\\
{\small P.O. Box 19395-5746, Tehran, Iran }\\[2pt]
{edwin.vandam@uvt.nl \quad romidi@cc.iut.ac.ir}}

\date{}

\maketitle \footnotetext[1] {This research is partially
carried out in the IPM-Isfahan Branch and in part supported
by a grant from IPM (No. 93050217). This version will appear in J. Alg Combinatorics, \url{https://doi.org/10.1007/s10801-017-0789-8}} \vspace*{-0.5cm}

\begin{abstract}
We generalize the concept of strong walk-regularity to directed graphs. We call a digraph
strongly $\ell$-walk-regular with $\ell >1$ if the number of walks of
length $\ell$ from a vertex to another vertex depends only on whether the first vertex is the same as, adjacent to, or not adjacent to the second vertex. This generalizes also the well-studied strongly regular digraphs and a problem posed by Hoffman.
Our main tools are eigenvalue methods. The case that the adjacency matrix is diagonalizable with only real eigenvalues resembles
the undirected case. We show that a digraph $\G$ with only real eigenvalues whose adjacency matrix is not diagonalizable has at most two values of $\ell$ for which $\G$ can be strongly $\ell$-walk-regular, and we also construct examples of such strongly walk-regular digraphs. We also consider digraphs with nonreal eigenvalues. We give such examples and characterize those digraphs $\G$ for which there are infinitely many $\ell$ for which $\G$ is strongly $\ell$-walk-regular. \\


\end{abstract}

\section{Introduction}
In \cite{vDO}, we introduced the concept of ``strongly walk-regular graphs'' as a generalization of
strongly regular graphs. Here we generalize this concept to directed graphs: a digraph is called {\em strongly $\ell$-walk-regular}
with $\ell >1$ if the number of walks of
length $\ell$ from a vertex to another vertex depends only on whether the first vertex is the same as, adjacent to, or not adjacent to the second vertex. In the undirected case, we used eigenvalues to characterize such graphs and we constructed several families of examples.
 Eigenvalue methods also play a crucial role in this paper, but the situation is more complex and interesting for two reasons. First of all, the adjacency matrix of a strongly walk-regular digraph need not be diagonalizable, and secondly, the eigenvalues can be nonreal.

The concept of strongly walk-regular digraphs generalizes the concept of strongly regular digraphs introduced by Duval \cite{Du} and is also related to a problem posed by Hoffman (unpublished).
He posed the problem of constructing digraphs with unique walks of length $3$, a problem that was generalized by Lam and Van Lint \cite{LvL} to arbitrary given length. Such digraphs are indeed special cases of strongly walk-regular graphs. Related work has also been done by
Bos\'{a}k \cite{Bosak} and Gimbert \cite{Gi}, who considered digraphs with unique walks of length in a given interval. Related is also the work of Comellas, Fiol, Gimbert, and Mitjana \cite{CFGM}, who introduced weakly distance-regular digraphs as those digraphs for which the number of walks of length $\ell$ from one vertex to another depends only on the distance between the vertices and on $\ell$. These generalize the (standard) distance-regular digraphs that were introduced by Damerell \cite{Damerell} and strongly regular digraphs.

In Section \ref{sec:pre} we shall give some background on eigenvalues of digraphs and strongly regular digraphs. In Section \ref{sec:3ev}, we will observe that in general, strong regularity of a digraph is a property that cannot be derived from the spectrum. This indicates behavior that is quite different from that of the undirected case. Indeed, we will give examples of non-strongly regular digraphs whose adjacency matrix is not diagonalizable, but still has the same spectrum as a particular strongly regular digraph. We will also use these examples in Section \ref{sec:realnondiagonalizable} to construct strongly walk-regular digraphs whose adjacency matrix is not diagonalizable. After making some basic observations in Section \ref{sec:basics}, we classify in Section \ref{sec:mu0} the strongly connected strongly $\ell$-walk-regular digraphs for which the number of walks of length $\ell$ from vertices to non-adjacent vertices is zero. In the remaining sections, we focus on the general case. In Section \ref{sec:regular}, we derive properties of the eigenvalues, which we use in Section \ref{sec:realeigenvalues}, which is on digraphs with real eigenvalues only. In Section \ref{sec:realdiagonalizable}, we focus on those digraphs for which the adjacency matrix is diagonalizable with real eigenvalues only. The results and examples in this case resemble those for undirected graphs as given in \cite{vDO}. In Section \ref{sec:realnondiagonalizable}, we show that a digraph $\G$ with only real eigenvalues and whose adjacency matrix is not diagonalizable has at most two values of $\ell$ for which $\G$ can be strongly $\ell$-walk-regular, and we also construct examples of such strongly walk-regular digraphs. In the final section, we focus on the digraphs with nonreal eigenvalues. We give examples and characterize those digraphs $\G$ for which there are infinitely many $\ell$ for which $\G$ is strongly $\ell$-walk-regular.

\section{Preliminaries}\label{sec:pre}

A {\it digraph} (or {\it directed graph}) $\G$ is an ordered pair $(V,E)$
consisting of a set $V$  of vertices and a set $E$ of ordered pairs of elements of $V$, called (directed) edges. We say that a vertex $u$ is adjacent to $v$
if the ordered pair $uv$ is an edge. In this case, we also call $v$ an outneighbor of $u$, and $u$ an inneighbor of $v$.
When $E$ contains both edges $uv$ and $vu$, we say that the edge $uv$ is bidirected. A digraph for which all edges are bidirected is considered the same as an {\it undirected graph}. A digraph having no multiple edges or loops (edges of the form $uu$) is called {\it simple}. All digraphs we consider in this paper are simple.

The {\it adjacency matrix} $A$ of a digraph $\G$ is the $n \times n$ matrix $(a_{uv})$ indexed by the vertices of $\G$,
with entries $a_{uv}=1$ if $u$ is adjacent to $v$, and $a_{uv}= 0$ otherwise. The all-ones matrix is denoted by $J$, or $J_n$ if we want to specify that its size is $n \times n$.
A digraph is called {\it regular} of degree $k$ if
$AJ = JA = kJ$, that is, if all vertices have both indegree and outdegree $k$. A {\it walk} (of length $\ell$) is a sequence of vertices $(u_0,u_1,\dots,u_{\ell})$, where $u_iu_{i+1}$ is an edge for $i=0,1,\dots,\ell-1$. The number of walks of length $\ell$ from $u$ to $v$ is given by $(A^{\ell})_{uv}$.  A digraph is {\it strongly connected} if there is a walk
from every vertex to every other vertex. The {\it reverse} of a digraph $\G$ is the digraph with adjacency matrix $A^{\top}$.

The line digraph of a graph $(V,E)$ has vertex set $E$. If $uv$ and $wz$ are both in $E$, then $uv$ is adjacent to $wz$ in the line digraph if $v=w$.

\subsection{The spectrum of a digraph}\label{sec:spectrum}

The {\it spectrum} of a digraph consists of the set of eigenvalues of its adjacency matrix together with their (algebraic) multiplicities.
Some basic results on the spectrum are the following (see also e.g. \cite{Brualdi}):
\begin{itemize}
\item[\rm(i)] By the Perron-Frobenius theorem, the maximum eigenvalue $\theta_0$ of a strongly connected
    digraph $\G$ is real, simple, and has a positive eigenvector ${\bf u}$. In particular, if $\G$ is $k$-regular
    then ${\bf u}=j$, where $j$ denotes the all-ones vector, and $\theta_0=k$;

\item[\rm(ii)] For a $k$-regular digraph $\G$, every eigenvector {\bf x} of an eigenvalue $\theta$ different from $k$ is orthogonal to the all-ones vector $j$. This follows from the equation $\theta j^{\top}{\bf x}=j^{\top}A{\bf x}=kj^{\top}{\bf x}$;

\item[\rm(iii)] If $\G$ is a strongly connected digraph with minimal polynomial having degree $d+1$, then the diameter of $\G$ is at most $d$;

\item[\rm(iv)] If $\G$ has $n$ vertices and $m$ edges, then the spectrum of the line digraph of $\G$ consists of the spectrum of $\G$ and $m-n$ extra eigenvalues $0$;

\item[\rm(v)] Hoffman and McAndrew \cite{HM} showed that for a digraph $\Gamma$ with adjacency matrix $A$, there exists a polynomial $h(x)\in
    \mathbb{Q}[x]$ such that
\begin{equation}\label{hoffmanpol}
h(A)=J
\end{equation}
if and only if $\Gamma$ is strongly connected and regular. If this is the case, then the unique polynomial $h(x)$ of least degree such that \eqref{hoffmanpol} is satisfied, the
    {\it Hoffman polynomial}, is $\frac n{\hbar(k)}\hbar(x)$, where $(x-k)\hbar(x)$ is the minimal polynomial of $\Gamma$
    and $n$ and $k$ are the number of vertices and the degree of $\Gamma$, respectively.
\end{itemize}

A useful consequence of the above is that if $\G$ is a strongly connected $k$-regular digraph with adjacency matrix $A$, and $p$ is a polynomial such that $p(A)=\alpha J$ for some $\alpha$, then the Hoffman polynomial divides $p$, and $p(\theta)=0$ for every eigenvalue $\theta$ different from $k$.

\subsection{Walk-regular and strongly regular digraphs}

Like in the undirected case, a digraph is called {\it walk-regular} if for every $\ell$, the number of (closed) walks of length $\ell$ from a vertex to itself is independent of the chosen vertex. This is equivalent to the property that $A^{\ell}$ has constant diagonal for every $\ell$.

The notion of a {\it strongly regular digraph} (or {\it directed strongly regular graph}) was introduced by
Duval \cite{Du} in 1988 as a generalization of strongly regular graphs to the directed
case. A strongly regular digraph is a regular digraph such that the number of walks of length two from one vertex to another
depends only on whether the first vertex is the same as, adjacent to, or not adjacent to the second vertex. In particular, a $k$-regular digraph on $n$ vertices with adjacency matrix $A$ is strongly regular with parameters $(n, k, t, \lambda, \mu)$ if $A^2= tI +\lambda A +\mu (J-I-A)$.  The case $t=k$ is the undirected case. On the other extreme, the case $t=0$, we have (doubly regular) tournaments, in which case $A+A^{\top}=J-I$. These two cases are typically excluded from the study of strongly regular digraphs. We therefore say that a strongly regular digraph is nonexceptional if $0<t<k$. Note that by the given definition, also the complete digraph (with adjacency matrix $J-I$) is strongly regular. We make this specific remark because in the undirected case, the complete graph is excluded from the definition of strong regularity. For more details, construction methods, and references, we refer to Brouwer's website \cite{Br}.

\subsection{Not a spectral characterization of strong regularity}\label{sec:3ev}

A connected regular undirected graph with three distinct eigenvalues is strongly regular. This does however not generalize to digraphs, as the next examples will show.
The smallest (nonexceptional) strongly regular digraph is on six vertices with spectrum $\{2^1,0^3,-1^2\}$; see \cite[Fig.~1]{Du}. There are however three other digraphs with this spectrum (which follows from checking all $2$-regular digraphs on six vertices that were generated by Brinkmann \cite[private communication]{Brinkmann}). Each of these three is strongly connected and regular, but the Hoffman polynomial is $x^2(x+1)$ (so it has $0$ as a multiple root), and therefore it is not strongly regular. The first of these is obtained from a directed $6$-cycle $\{12,23,34,45,56,61\}$ by adding bidirected edges $13,25,46$. The second digraph has directed edges $12,23,31,14,45,51$ and bidirected edges $24,36,56$, whereas the third digraph is the reverse of the second. The latter two digraphs are not even walk-regular (because $A^2$ does not have constant diagonal), but the first one is. Thus, we have the following.

\begin{observation} Strong regularity of digraphs can in general not be recognized from the spectrum.
\end{observation}

Note that because the Hoffman polynomials of the above examples have $0$ as a multiple root, their adjacency matrices are not diagonalizable. Godsil, Hobart, and Martin \cite{GHM} showed that for nonexceptional strongly regular digraphs, the adjacency matrix is diagonalizable. On the other hand, a regular digraph with three distinct eigenvalues whose adjacency matrix is diagonalizable must be strongly regular. Later we will also see examples of strongly walk-regular digraphs for which the adjacency matrix is not diagonalizable. The fact that this can happen is one of the interesting features in the generalization of undirected graphs to directed graphs.

\section{The basics}\label{sec:basics}

Like in the undirected case, we call a digraph $\G$ a {\em strongly $\ell$-walk-regular digraph with parameters} $(\lal,\mul,\nul)$, for $\ell>1$, if the number of walks of length $\ell$ from a vertex to an adjacent vertex equals $\lal$, from a vertex to a non-adjacent vertex equals $\mul$, and from a vertex to itself equals $\nul$. So indeed, every strongly regular
digraph with parameters $(n, k, t, \lambda, \mu)$ is a strongly $2$-walk-regular digraph with parameters
$(\lambda,\mu,t)$. In particular, the empty and complete digraphs
give examples. Indeed, these digraphs are clearly strongly $\ell$-walk-regular for every $\ell$. It is also clear that if $\G$ is an undirected digraph (i.e., its adjacency matrix is symmetric), then it is strongly $\ell$-walk-regular as an undirected graph in the sense of \cite{vDO} if and only if it is strongly $\ell$-walk-regular as a digraph in the above sense.

We will first make some basic observations. All of these are similar as in the undirected case, so we omit the (elementary) proofs.

Let $A$ be the adjacency matrix of $\G$. Then $\G$ is a strongly
$\ell$-walk-regular digraph if and only if $A^{\ell}$ is in the span
of $A,I$, and $J$.

\begin{lemma}\label{lem:A} Let $\ell>1$, and let $\G$ be a digraph with adjacency matrix $A$.
Then $\G$ is a strongly $\ell$-walk-regular digraph with parameters
$(\lal,\mul,\nul)$ if and only if $A^{\ell}+(\mul-\lal)A+(\mul-\nul)I=\mul J$.
\end{lemma}

Now it is clear that a strongly regular digraph is strongly $\ell$-walk-regular for more values of $\ell$ than just $\ell=2$, because its adjacency algebra (that is, the algebra spanned by all powers of $A$) equals $\langle A,I,J \rangle$.

\begin{proposition}\label{prop:srg} Let $\G$ be a strongly regular digraph.
Then $\G$ is a strongly $\ell$-walk-regular digraph with parameters $(\lal,\mul,\nul)$
for every $\ell>1$ and some $\lal,\mul,$ and $\nul$.
\end{proposition}

It is also clear from Lemma \ref{lem:A} that the reverse of a strongly walk-regular graph is strongly walk-regular, with the same parameters.

By Hoffman and McAndrew's characterization of strongly connected regular digraphs \cite{HM} (see Section \ref{sec:spectrum}), we have the following.

\begin{lemma} \label{strongly connected} Let $\ell>1$, and let $\G$ be a
strongly $\ell$-walk-regular digraph with parameters $(\lal,\mul,\nul)$ where
$\mul>0$. Then $\G$ is regular and strongly connected.
\end{lemma}

Also if $\mul=0$, the digraph can be regular and strongly connected. For example, the directed
cycle of size $\ell$ is strongly $\ell$-walk-regular
with parameters $(0,0,1)$ and strongly $(\ell+1)$-walk-regular
with parameters $(1,0,0)$. We will look further into the case $\mul=0$ in Section \ref{sec:mu0}, and focus on the case of strongly connected regular digraphs in the later sections.

\section{Graphs with $\mu_{\ell}=0$}\label{sec:mu0}

In this section, we shall classify the strongly $\ell$-walk-regular digraphs with $\mul=0$ that are strongly connected.

\begin{example}\label{ex:mu01} Let $\G$ be a directed $g$-cycle, with adjacency matrix $A$. It is clear that $\G$ has $x^{g-1}+\cdots+x+1$ as its Hoffman polynomial and $x^{g}-1$ as its minimal polynomial. The eigenvalues of $\G$ are the complex $g^{\text{th}}$-roots of unity. From the equation $A^{g}=I$, it follows that $\G$ is strongly $\ell$-walk-regular, with $\mul=0$, for $\ell \equiv 0$ and $1 \text{ (mod }g)$.
\end{example}

By using the so-called coclique extension, we can construct more examples as follows.

\begin{example}\label{ex:mu02} Let $\G$ be a coclique extension of a directed $g$-cycle with $g \geq 2$, that is, a digraph with vertex set $V=\cup_{i=1}^{g} V_i$ and edge set $E=\cup_{i=1}^{g-1} (V_i\times V_{i+1})\cup (V_{g}\times V_{1})$. Let $A$ be the adjacency matrix of $\G$, then $A^{g+1}=\lambda_{g+1}A$, where $\lambda_{g+1}=\Pi_{i=1}^{g}|V_i|$. If we require that $|V_i|>1$ for at least one $i$, then $\G$ has diameter $g$, which implies that $x^{g+1}-\lambda_{g+1}x$ is the minimal polynomial of $\G$. It follows that $\G$ is strongly $\ell$-walk-regular, with $\mul=0$, for $\ell \equiv 1 \text{ (mod }g)$. Note that the case $g=2$ gives a (undirected) complete bipartite graph.
\end{example}

The properties of these two examples are very typical for the strongly walk-regular digraphs that have nonreal eigenvalues, as we shall see in Section \ref{sec:nonreal}. We now first mention some exceptional nonregular examples.

\begin{example}\label{ex:mu0exceptional} Let $\ell \geq 3$, and consider the directed $\ell$-cycle on vertex set $\mathbb{Z}_{\ell}$, where a vertex $u$ is adjacent to a vertex $v$ if $v=u+1$. To this digraph, we add the edge $02$. The obtained digraph has diameter $\ell-1$, and is strongly $\ell$-walk-regular with parameters $(\lal,\mul,\nul)=(1,0,1)$. It follows that its minimal polynomial is $x^{\ell}-x-1$. If on top, we also add the edge $13$, then the obtained digraph is strongly $\ell$-walk-regular with parameters $(\lal,\mul,\nul)=(2,0,1)$, and its minimal polynomial is $x^{\ell}-2x-1$. Note that for $\ell=3$, the latter digraph is the complete graph $K_3$ minus an edge.
\end{example}

We will show next that the given examples are the only strongly connected strongly $\ell$-walk-regular digraphs with $\mul=0$. In order to do this, we will use the following lemma on shortest cycles in such digraphs. We recall that the girth $g$ of a digraph is the length of the shortest directed cycle. Note that $g=2$ if there is a bidirected edge.

\begin{lemma}\label{lem:girth} Let $\G$ be a strongly connected strongly $\ell$-walk-regular digraph with $\mul=0$ and girth $g$. Then $\ell \equiv 0$ or $1 \textup{ (mod }g)$. Moreover, if $\G$ is not a directed cycle, then $\lal>0$ and for every directed cycle $C$ of length $g$ and every vertex $z$ not on $C$ there are vertices $u$ and $v$ on $C$ such that $uz$ and $zv$ are edges.
\end{lemma}

\pf
Let $C$ be a cycle of shortest length $g$, and let $w$ be a vertex on $C$. It is clear that if one starts walking on the cycle from $w$, then after $\ell$ steps one should end up in either $w$ itself or the unique outneighbor of $w$ on $C$, because $\mul=0$. Thus $\ell \equiv 0$ or $1 \text{ (mod }g)$.

Now suppose that $\G$ is not a directed cycle, and let $z$ be a vertex that is not on $C$. Because $\G$ is strongly connected and the diameter of $\G$ is less than $\ell$ (because the minimal polynomial is a divisor of $x^{\ell}-\lal x- \nul$), it follows that there is a vertex on $C$ which is at distance less than $\ell$ from $z$. But then one can extend a walk from $z$ to that vertex to a walk of length $\ell$ from $z$ to a vertex $v$ on $C$. Because $\mul=0$, it follows that $zv$ is an edge and that $\lal>0$. Similarly there is a vertex $u$ on $C$ such that $uz$ is an edge.
\epf

\begin{theorem}\label{prop:mu0} Let $\G$ be a strongly connected strongly $\ell$-walk-regular digraph with $\mul=0$ that is not complete. Then $\G$ is one of the digraphs of Examples \ref{ex:mu01}-\ref{ex:mu0exceptional}.
\end{theorem}

\pf Note again that $\G$ has diameter less than $\ell$, and so we may assume that $\ell>2$. In the following, let us denote by $a \rightarrow_i b$ that there is a walk of length $i$ from $a$ to $b$; we omit the subscript for $i=1$. Let $g$ be the girth of $\G$ and consider a cycle $C$ of length $g$. By Lemma \ref{lem:girth}, we know that $\ell \equiv 0$ or $1 \textup{ (mod }g)$. Further, we may assume that $\G$ is not a $g$-cycle (Example \ref{ex:mu01}). Let $z$ be a vertex not on $C$, and let $u$ and $v$ be vertices on $C$ such that $uz$ and $zv$ are edges (see Lemma \ref{lem:girth}).

First, assume that $\ell \equiv 0 \textup{ (mod }g)$. Let $w$ be the outneighbor of $u$ on $C$. By walking $\ell-1$ steps from $w$ on the cycle, one ends up in $u$. So $w \rightarrow_{\ell-1} u \rightarrow z$, and so there is a walk of length $\ell$ from $w$ to $z$, hence $wz$ is also an edge. Inductively, it follows that every vertex on $C$ is adjacent to $z$. Similarly, $z$ is adjacent to every vertex on $C$. This also implies that $g=2$, and it easily follows that $\G$ is a complete graph. Thus, we may assume below that $\ell \equiv 1 \textup{ (mod }g)$.

Let us now assume that $\nul=0$. We again consider the above vertices $z$, $u$, and $v$. Because $C$ is a shortest cycle, it follows that the distance from $u$ to $v$ on the cycle is at most two. If $uv$ is an edge, then $v \rightarrow_{\ell-2} u \rightarrow z \rightarrow v$, so there is a closed walk of length $\ell$, which contradicts the assumption that $\nul=0$. Thus, $u$ and $v$ are at distance two, or $u=v$ in the (degenerate) case that $g=2$. In any case, given $C$ and $z$, the vertices $u$ and $v$ are the unique vertices on $C$ such that $uz$ and $zv$ are edges. Now let $v_i$, for $i=1,\dots,g$ be the consecutive vertices of $C$,  so $v_{i-1}v_{i}$ is an edge for $i=1,\dots,g$, where we let $v_0=v_g$. We now define the set of vertices $V_i=\{v : v_{i-1}v \text{ is an edge}\}$ for $i=1,\dots,g$. It is easy to show now that there are no edges within each of the sets $V_i$, that each vertex in $V_{i}$ is adjacent to each vertex in $V_{i+1}$, and finally to draw the conclusion that $\G$ must be a coclique extension of a $g$-cycle (Example \ref{ex:mu02}).

Next, we assume that $\nul>0$ (and recall that we assumed that $\ell \equiv 1 \textup{ (mod }g)$), and hence that there are closed walks of length $\ell$. Any such closed walk must contain a vertex not on $C$. In particular, it follows that there is an edge $uz$, with $u$ on $C$ and $z$ not on $C$, that is contained in a closed walk of length $\ell$. If $u$ has outneighbor $v$ on $C$, then it follows that $z \rightarrow_{\ell-1} u \rightarrow v$, so then $zv$ is an edge.
It also follows that every edge of $C$, except possibly $uv$ is contained in a closed walk of length $\ell$. Suppose now that besides $z$, there is another vertex $z'$ not on $C$, and let $v'$ be an outneighbor of $z'$ on $C$. If $v' \neq v$, and $v''$ is the inneighbor of $v'$ on $C$, then $v''v'$ is contained in a walk of length $\ell$, so $z' \rightarrow v' \rightarrow_{\ell-1} v''$, and so $z'v''$ is also an edge. Inductively, it follows that $z'v$ is an edge. By a similar argument (or applying the same argument to the reverse digraph), it follows that $uz'$ is an edge. Now $z \rightarrow v \rightarrow_{\ell-2} u \rightarrow z'$, so $zz'$ is an edge. By the same argument, this edge is bidirected, and so $g=2$. It is now easy to show that the digraph is complete, which is a contradiction. Hence there is only one vertex ($z$) not on $C$. If $z$ has only one outneighbor ($v$) and one inneighbor ($u$), then $\G$ is the first digraph of Example \ref{ex:mu0exceptional}. Because $C$ is a cycle of shortest length, $z$ can have either one more inneighbor (the inneighbor of $u$ on $C$) or one more outneighbor (the outneighbor of $v$ on $C$), but not both. In both cases, $\G$ is the second digraph of Example \ref{ex:mu0exceptional}.
\epf

Note that for the digraphs of Example \ref{ex:mu0exceptional}, the above proof indicates that $\ell \equiv 1 \textup{ (mod }g)$. It is however not hard to show that $\ell=g+1$ (by constructing walks of length $tg+1$ from a vertex to a non-adjacent vertex for $t \geq 2$), i.e., there is only one $\ell$ for which these digraphs are strongly $\ell$-walk-regular.

The given examples with $\mul=0$ and $\nul=0$ can be used to construct also some examples that are weakly connected. For example, take a directed cycle, and add a few vertices that only have one outneighbor (and no inneighbors), and this outneighbor is on the cycle. Even more degenerate examples can be constructed: if all walks eventually end in vertices without outneighbor, then $A^{\ell}=O$ for $\ell$ large enough.

\section{Strongly connected regular digraphs}\label{sec:regular}

In order to study the case that $\mul>0$, from now on we consider strongly connected and regular digraphs; see Lemma \ref{strongly connected}. We denote the set of all diagonalizable digraphs (that is, digraphs with diagonalizable adjacency matrix) by $\mathcal{D}$. For these digraphs, all roots of the minimal polynomial are simple. By $\mathcal{D}_{\theta}$  we denote the set of digraphs whose minimal polynomial has all but one root simple, and the nonsimple root is $\theta$ and it has multiplicity $2$. For example, the three digraphs in Section \ref{sec:3ev} that are cospectral to a strongly regular digraph, are in $\mathcal{D}_{0}$. The following result shows that if $\G$ is
strongly walk-regular, then either $\G\in \mathcal{D}$ or $\G\in \mathcal{D}_{\theta}$ for some eigenvalue $\theta$ of $\G$. Note that from the observations in Section \ref{sec:spectrum} and Lemma \ref{lem:A}, it follows that the Hoffman polynomial divides the polynomial $x^{\ell}+(\mul-\lal)x+\mul-\nul$, and hence each eigenvalue different from the degree $k$ is a root of this polynomial.

\begin{proposition}\label{multiplicities}
Let $\G$ be a strongly connected $k$-regular digraph. If  $\G$ is strongly $\ell$-walk-regular
with parameters $(\lal,\mul,\nul)$ where $\ell>1$, then either $\G\in \mathcal{D}$ or $\G\in\mathcal{D}_{\theta}$ for $\theta=\frac{-\ell(\mul-\nul)}{(\ell-1)(\mul-\lal)}$. Moreover, if $\G\in\mathcal{D}_{\theta}$, then $\theta$ is a nonzero integer number different from $k$, and $(\frac{\mul-\nul}{\ell-1})^{\ell-1}=(\frac{\mul-\lal}{-\ell})^{\ell}$.
\end{proposition}

\pf Suppose that $\G$ is not diagonalizable and let $\theta$ be a nonsimple root of the minimal polynomial. Then $\theta \neq k$ and hence it is also a nonsimple root of the Hoffman polynomial, and hence of $p(x)=x^{\ell}+(\mul-\lal)x+\mul-\nul$. If the multiplicity of $\theta$ in these polynomials is larger than two, then clearly it is also a root of the derivatives $p'(x)=\ell x^{\ell-1}+\mul-\lal$ and
$p''(x)=\ell (\ell-1) x^{\ell-2}$. This implies that $\theta=0$ and hence that $\mul=\lal=\nul$. Therefore $A^{\ell}=\mul J$ and so every eigenvalue different from $k$ is $0$. This however contradicts the fact that $A$ has trace $0$, so $\theta$ must have multiplicity two in $p(x)$, and hence in the minimal polynomial.

If $\theta$ is indeed a root of $p(x)$ with multiplicity two, then it is also a root of $p'(x)$, and similar as before it follows that $\theta \neq 0$. By combining the equations $p(\theta)=0$ and $p'(\theta)=0$, it follows that
\begin{equation}\label{eq:mullal}
\mul-\lal=-\ell \theta^{\ell-1}\textrm{~~and~~}\mul-\nul=(\ell-1) \theta^{\ell}
\end{equation} and hence that $\theta=\frac{-\ell(\mul-\nul)}{(\ell-1)(\mul-\lal)}$. If this is indeed a root of both polynomials, then
$(\frac{\mul-\nul}{\ell-1})^{\ell-1}=(\frac{\mul-\lal}{-\ell})^{\ell}$. Moreover, if this is the case, then $\theta$ is a rational eigenvalue, and so it must be integer (being a root of a monic polynomial with integer coefficients).\epf

In Section \ref{sec:realnondiagonalizable}, we will show that the case of nondiagonalizable strongly walk-regular digraphs really occurs by constructing some examples.

\begin{proposition}\label{eigenvalues}
Let $\ell>1$. A strongly connected $k$-regular digraph $\G$ on $n$ vertices is strongly $\ell$-walk-regular
with parameters $(\lal,\mul,\nul)$ if and only if all of the following conditions hold:
\begin{itemize}
\item[\rm(i)] Either $\G\in \mathcal{D}$ or $(\frac{\mul-\nul}{\ell-1})^{\ell-1}=(\frac{\mul-\lal}{-\ell})^{\ell}$ and $\G\in \mathcal{D}_{\theta}$ for
$\theta=\frac{-\ell(\mul-\nul)}{(\ell-1)(\mul-\lal)};$
\item[\rm(ii)] All eigenvalues besides $k$ are roots of the
equation
\begin{equation*}x^{\ell}+(\mul-\lal)x+\mul-\nul=0;\end{equation*}
\item[\rm(iii)]
And \begin{equation*}k^{\ell}+(\mul-\lal)k+\mul-\nul=\mul n.\end{equation*}
\end{itemize}
\end{proposition}

\pf If $\G$ is a strongly $\ell$-walk-regular digraph with parameters
$(\lal,\mul,\nul)$, then condition $\rm(i)$ holds by Proposition \ref{multiplicities} and $A^{\ell}+(\mul-\lal)A+(\mul-\nul)I=\mul J$ by Lemma \ref{lem:A}. We already observed in Section \ref{sec:spectrum} that this implies $\rm(ii)$. Condition $\rm(iii)$ follows from multiplying the
above matrix equation with the all-ones vector.

Now assume that $\rm(i), \rm(ii)$ and $\rm(iii)$ hold. If $\G$ is diagonalizable, then each eigenvalue besides $k$ is a simple root of the Hoffman polynomial, and hence the Hoffman polynomial divides the polynomial $p(x)=x^{\ell}+(\mul-\lal)x+\mul-\nul$. If $(\frac{\mul-\nul}{\ell-1})^{\ell-1}=(\frac{\mul-\lal}{-\ell})^{\ell}$ and $\G\in \mathcal{D}_{\theta}$ for
$\theta=\frac{-\ell(\mul-\nul)}{(\ell-1)(\mul-\lal)}$, then $\theta$ is a root of both $p(x)$ and $p'(x)$, and so it is a root of multiplicity at least two of $p(x)$. Because all other eigenvalues besides $k$ and $\theta$ are simple roots of the Hoffman polynomial, and $\theta$ is a root of multiplicity two in the Hoffman polynomial, also in this case the Hoffman polynomial $h(x)$ divides $p(x)$.

Now let $g(x)=p(x)/h(x)$. Because $g(x)$ is a polynomial, it follows that $p(A)=g(A)h(A)=g(A)J=g(k)J$. By multiplying with the all-ones vector and condition $\rm(iii)$, this implies that $\mul n= p(k)=g(k)n$, and hence we obtain that $p(A)=\mul J$. Thus $\G$ is strongly $\ell$-walk-regular with parameters $(\lal,\mul,\nul)$.\epf

Like in the undirected case, we obtain the following result from the above proof.

\begin{corollary}\label{cor:eigenvalues} Let $\ell>1$. A strongly connected regular digraph is strongly
$\ell$-walk-regular if and only if its Hoffman polynomial divides the
polynomial $x^{\ell}+ex+f$ for some integers $e$ and $f$.
\end{corollary}

Because we can bound the number of real roots of $x^{\ell}+ex+f$, this has consequences for the
number of real roots of the Hoffman polynomial, and hence for the number of distinct real eigenvalues of a strongly walk-regular digraph. The bound is as follows.

\begin{lemma} \label{lem:descartes} {\em \cite{vDO}} Let $\ell>1$, and let $p(x)=x^{\ell}+ex+f$ for some real $e$ and $f$. Then $p$ has at most three real roots. If $\ell$ is even, then $p$ has at
most two real roots.
\end{lemma}

Note that the counted number of real roots includes multiplicities of these roots. We thus obtain the following from Corollary \ref{cor:eigenvalues}.

\begin{proposition}\label{prop:3or4} Let $\G$ be strongly connected, regular, and
strongly $\ell$-walk-regular with $\ell>1$. Then the Hoffman polynomial of $\G$ has at most three
real roots. Moreover, if $\ell$ is even, then it has at most two
real roots.
\end{proposition}

We will now give three more results that follow from Proposition \ref{eigenvalues}. We will use these in the next sections.

\begin{proposition}\label{charactrization1} Let $\G$ be a strongly connected $k$-regular digraph on $n$ vertices with at least three distinct eigenvalues, and let $\ell>1$. Then $\G$ is strongly $\ell$-walk-regular with parameters $(\lal,\mul,\nul)$ if and
only if all of the following conditions hold:
\begin{itemize}
\item[\rm(i)] Either $\G\in \mathcal{D}$ or $(\frac{\mul-\nul}{\ell-1})^{\ell-1}=(\frac{\mul-\lal}{-\ell})^{\ell}$ and $\G\in \mathcal{D}_{\theta}$ for
$\theta=\frac{-\ell(\mul-\nul)}{(\ell-1)(\mul-\lal)};$
\item[\rm(ii)] For every two distinct eigenvalues
$\theta_1,\theta_2\neq k$,
\begin{equation}\label{eq:lswr1}
\mul-\lal=\frac{\theta_2^{\ell}-\theta_1^{\ell}}{\theta_1-\theta_2} \textrm{~~and~~} \mul-\nul=\frac{\theta_2\theta_1^{\ell}-\theta_1\theta_2^{\ell}}{\theta_1-\theta_2};
\end{equation}
\item[\rm(iii)]
And \begin{equation*}k^{\ell}+(\mul-\lal)k+\mul-\nul=\mul n.\end{equation*}
\end{itemize}
\end{proposition}

\pf Let $\theta_1,\theta_2\neq k$ be two distinct eigenvalues of $\G$. Then it is straightforward to show that \eqref{eq:lswr1} holds if and only if $\theta_1$ and $\theta_2$ are roots of the equation $x^{\ell}+(\mul-\lal)x+\mul-\nul=0$. This implies that condition (ii) of Proposition \ref{eigenvalues} is equivalent to the property that \eqref{eq:lswr1} holds for every two distinct eigenvalues $\theta_1,\theta_2\neq k$, which is all we have to show.\epf

Note that the restriction on $\G$ having at least three distinct eigenvalues is not really necessary. It is not so hard to see that a strongly connected digraph with (at most) two distinct eigenvalues must be a complete digraph (i.e., $A=J-I$), and this satisfies the conditions of the proposition.

\begin{corollary}\label{cor:three} Let $\G$ be a strongly connected $k$-regular digraph with at least four distinct eigenvalues, and let $\ell>1$. If $\G$ is strongly $\ell$-walk-regular, then
\begin{equation}\label{eq:lswr}
(\theta_2-\theta_3)\theta_1^{\ell}+(\theta_3-\theta_1)\theta_2^{\ell}+(\theta_1-\theta_2)\theta_3^{\ell}=0
\end{equation}
for every three distinct eigenvalues $\theta_1,\theta_2,\theta_3\neq k$.
\end{corollary}

\pf From \eqref{eq:lswr1}, it follows that $\frac{\theta_2^{\ell}-\theta_1^{\ell}}{\theta_1-\theta_2}=\frac{\theta_3^{\ell}-\theta_1^{\ell}}{\theta_1-\theta_3}$, and by working this out, \eqref{eq:lswr} follows. \epf

\begin{lemma}\label{lemma:nondiag} Let $\G$ be a strongly connected $k$-regular digraph that is strongly $\ell$-walk-regular with $\ell>1$. If $\G \in \mathcal{D}_{\theta}$, then
\begin{equation}\label{eq:eta}
\ell \theta^{\ell-1}(\theta-\eta)=\theta^{\ell}-\eta^{\ell}
\end{equation}
for every eigenvalue $\eta \neq k$.
\end{lemma}

\pf From Proposition \ref{eigenvalues}, it follows that $\eta^{\ell}+(\mul-\lal)\eta+\mul-\nul=0$. By using \eqref{eq:mullal}, it now follows that $\ell \theta^{\ell-1}(\theta-\eta)=\theta^{\ell}-\eta^{\ell}$.
\epf

\section{Digraphs with all eigenvalues real}\label{sec:realeigenvalues}
In this section, we will consider the digraphs that have real eigenvalues only. The case that $\ell$ is even is then easy.

\begin{proposition}\label{prop:diagrealeven} Let $\G$ be a strongly connected regular digraph with all eigenvalues real, and let $\ell$ be even. Then $\G$ is strongly $\ell$-walk-regular if and only if $\G$ is strongly regular.
\end{proposition}

\pf Assume that $\G$ strongly $\ell$-walk-regular. Because $\ell$ is even and all eigenvalues are real, it follows from Proposition \ref{prop:3or4} that the Hoffman polynomial of $\G$ has at most two roots (including multiplicities), and hence $\G$ is strongly regular. On the other hand, we already observed in Proposition \ref{prop:srg} that a strongly regular digraph is strongly $\ell$-walk-regular for every $\ell$. \epf

For $\ell$ odd, we will distinguish between diagonalizable digraphs and the others.

\subsection{Diagonalizable digraphs}\label{sec:realdiagonalizable}

Let $\G$ be diagonalizable with all eigenvalues real, and suppose that $\G$ is strongly connected, $k$-regular, and strongly $\ell$-walk-regular, but not strongly regular. Then $\ell$ is odd and it follows by Proposition \ref{prop:3or4} that $\G$ has four distinct eigenvalues $k>\theta_1>\theta_2>\theta_3$. The theory that was developed for strongly walk-regular undirected graphs with four eigenvalues in \cite[Section 4]{vDO} can almost literally be extended to this case. In particular, we obtain the following results.

\begin{proposition}
Let $\G$ be a strongly connected $k$-regular diagonalizable digraph with all eigenvalues real. If $\G$ is strongly $\ell$-walk-regular but not strongly regular, then $\ell$ is odd, $\G$ has four distinct real eigenvalues $k>\theta_1>\theta_2>\theta_3$, and
\begin{itemize}
\item[\rm(i)] $\G$ is strongly $3$-walk-regular if and only if $\theta_1+\theta_2+\theta_3=0$;
\item[\rm(ii)] If $\theta_2=0$ and
$\theta_3=-\theta_1$, then $\G$ is strongly $\ell$-walk-regular for every odd
$\ell$, and
\item[\rm(iii)]
If $\theta_2 \neq 0$ or $\theta_3 \neq
-\theta_1$, then there is at most one $\ell>1$ such that $\G$ is strongly
$\ell$-walk-regular.
\end{itemize}
\end{proposition}

\pf The arguments for $\rm(i)$, $\rm(ii)$, and $\rm(iii)$ are similar as those for \cite[Prop.~4.1]{vDO},  \cite[Prop.~4.2]{vDO}, and \cite[Thm.~4.4]{vDO}, respectively.
\epf

As in the undirected case, we can construct examples of strongly walk-regular digraphs by product constructions. Other examples are obtained by using line digraphs. The undirected version of this, the usual line graph, did not provide examples in the undirected case.

\begin{example}\label{ex:1} Consider a strongly regular digraph $\G$ with parameters $(n, k, t, \lambda, \mu)$, for which $\lambda=\mu \neq t$. An infinite family of such digraphs was constructed by J{\o}rgensen \cite{J} (see also \cite[Construction T4]{Br} and Example \ref{ex:4}). If $A$ is the adjacency matrix of such a digraph, and $q>1$, then the digraph with adjacency matrix $A \otimes J_q$ is diagonalizable with four distinct eigenvalues $qk$, $\pm q \sqrt{t-\mu}$, and $0$. So it is strongly $\ell$-walk-regular for every odd $\ell$.

Also the line digraph of $\G$ is diagonalizable with four distinct eigenvalues: $k$, $\pm \sqrt{t-\mu}$, and $0$, because the eigenvalues of a digraph and its line digraph only differ in the eigenvalues $0$, see Section \ref{sec:spectrum}. So also the line digraph of $\G$ is strongly $\ell$-walk-regular for every odd $\ell$. We remark finally that one could even combine the product construction and the line digraph construction (repeatedly).
\end{example}

\begin{example}\label{ex:2} Consider a strongly regular digraph $\G$ with parameters $(n, k, t, \lambda, \mu)$, for which $t>\mu=\lambda+1$ (see \cite[Construction M4]{Br} for an infinite family). Denote its eigenvalues by $k$, $r$, and $s$, then $r+s=-1$ and $r,s \neq 0$. Let $A$ be the adjacency matrix of such a digraph, then the digraph with adjacency matrix $J_{3n}-I-A \otimes J_3$ is diagonalizable with four distinct eigenvalues $3n-1-3k,3r+2,-1,-1-3r$. So it is a strongly $\ell$-walk-regular digraph only for $\ell=3$. Also here variations are possible by using line digraphs, for example by first taking the line digraph of $\G$, with adjacency matrix $B$, and then construct the digraph with adjacency matrix $J_{3nk}-I-B \otimes J_3$.
\end{example}

It would be interesting, as in the undirected case, to also find examples that are strongly $\ell$-walk-regular for precisely one $\ell$ with $\ell>3$ (and $\mul>0$; see Example \ref{ex:mu0exceptional} for such digraphs with $\mul=0$).

\subsection{Nondiagonalizable digraphs}\label{sec:realnondiagonalizable}

Using the earlier examples from Section \ref{sec:3ev} and the Kronecker product again, we can construct examples of nondiagonalizable strongly walk-regular digraphs.

\begin{example}\label{ex:3} Let $A$ be the adjacency matrix of a regular digraph on six vertices with Hoffman polynomial $x^2(x+1)$ (see Section \ref{sec:3ev} for three such examples). Then the digraph with adjacency matrix $J_{18}-I-A\otimes J_3$ has Hoffman polynomial $x^3-3x-2$. Indeed, it has distinct eigenvalues $11$, $2$, and $-1$, with $-1$ having multiplicity two in the minimal polynomial. So it is a strongly $3$-walk-regular digraph with three eigenvalues that is in $\mathcal{D}_{-1}$. It has the same spectrum as a strongly regular digraph that is constructed in the same way from the strongly regular digraph on six vertices that was mentioned in Section \ref{sec:3ev}. It can be shown that the nondiagonalizable digraph is strongly $\ell$-walk-regular only for $\ell=3$, following the approach of the proof of the next result.
\end{example}

\begin{proposition}\label{prop:realnondiag} Let $\G$ be a nondiagonalizable strongly connected regular digraph with all eigenvalues real. Then $\G$ is strongly $\ell$-walk-regular for at most two values of $\ell$, and these $\ell$ are odd.
\end{proposition}

\pf Assume that $\G \in \mathcal{D}_{\theta}$ is strongly $\ell$-walk-regular, and let $\eta \neq \theta,k$ be another real eigenvalue. Suppose that $\G$ is also strongly $m$-walk-regular, then $\ell \theta^{\ell-1}(\theta-\eta)=\theta^{\ell}-\eta^{\ell}$
and $m \theta^{m-1}(\theta-\eta)=\theta^{m}-\eta^{m}$ according to Lemma \ref{lemma:nondiag}. By combining these two equations, it follows that
\begin{equation*}
\ell(1-(\eta/\theta)^m)=m(1-(\eta/\theta)^{\ell}).
\end{equation*}
Now let $\xi =|\eta/\theta|$, $\epsilon= \sgn (\eta/\theta)$, and $f(x)=\ell(1-\epsilon \xi^x)-x(1-\epsilon \xi^{\ell})$. Then $f(m)=0$ and $f(\ell)=0$ because $m$ and $\ell$ are odd (by Proposition \ref{prop:diagrealeven}). We now aim to show that the equation $f(x)=0$ has at most two solutions, because that would finish the proof. Indeed, this follows from the fact that $f''(x)$ is always positive, or always negative, depending on $\epsilon$, unless $\eta=-\theta$. In this exceptional case however, $f(x)=2\ell-2x$, so it has only one root.
\epf

We would not be surprised if it can be shown that a nondiagonalizable strongly connected digraph with all eigenvalues real can be strongly $\ell$-walk-regular for at most one $\ell$. For the digraph of Example \ref{ex:3}, we find that the function $f$ from the above proof satisfies $f(x)=3(1+2^x)-9x$, and it is easy to show that this has only one integral root.

\section{Digraphs with nonreal eigenvalues}\label{sec:nonreal}

In this section, we consider strongly walk-regular digraphs with nonreal eigenvalues. Doubly regular tournaments, or equivalently strongly regular digraphs with $t=0$ are examples of these: If $n$ is the number of vertices, then the eigenvalues are $k=\frac12(n-1)$ and $-\frac12 \pm \frac12 \sqrt{-n}$ (this follows for example from the proof of \cite[Thm.~2.2]{Du}). Clearly also the directed cycle (see Example \ref{ex:mu01}) and the digraphs of Example \ref{ex:mu02} have nonreal eigenvalues.

Hoffman (unpublished) posed the problem of constructing digraphs with unique walks of length $3$. Lam and Van Lint \cite{LvL} generalized this by considering directed graphs with unique walks of fixed length ($m$ say), that is, with adjacency matrix $A$ satisfying $A^{m}=J-I$. This is a very particular case of our strongly walk-regular digraphs (with $\mu_m=\lambda_m=1$ and $\nu_m=0$), whose eigenvalues different from $k$ satisfy the equation $x^{m}=-1$. In particular, they showed that there are no such digraphs for even $m$, and constructed $k$-regular digraphs for every odd $m$ on $k^{m}+1$ vertices. In order to generalize their example, we use the following lemma.

\begin{lemma}\label{lemmaW} Let $\ell \geq 1$, $k \geq 2$, and $W=\{\sum_{i=0}^{\ell-1} a_i(-k)^i:a_0,a_1,\dots,a_{\ell-1} \in \{1,2,\dots,k\}\}$. Then
$W=\{1,2,\dots,k^{\ell}\}$ if $\ell$ is odd, and $W=\{0,-1,-2,\dots,-k^{\ell}+1\}$ if $\ell$ is even.
\end{lemma}

\pf This follows easily by induction on $\ell$.
\epf

Now we can generalize the construction of Lam and Van Lint \cite{LvL} to include even $m$.

\begin{example}\label{ex:4} Let $k \geq 2$, $m \geq 2$, and $n=k^{m}-(-1)^m$. The vertex set of $\G$ is $\mathbb{Z}_n$, and a vertex $u$ is adjacent to a vertex $v$ if $ku+v \in \{1,2,\dots,k\}$. This is an example of a so-called $(-k)$-circulant graph (see \cite{Lam}). Note that $n$ has a factor $k+1$, which implies that there are no loops, i.e., no vertex is adjacent to itself.

Lam and Van Lint \cite{LvL} showed that for $m=3$ and  $k=2$, this digraph $\G$ on nine vertices (which is depicted in their Fig.~3) is in fact the only $2$-regular digraph satisfying the equation $A^3=J-I$. Its eigenvalues are $2$, $-1$ with multiplicity $4$, and $\frac12 \pm \frac12 \sqrt{-3}$ with multiplicity $2$ each. It follows that $\G$ is strongly $\ell$-walk-regular for $\ell \equiv 0$ and $1 \text{ (mod 3)}$.

In general, observe (by induction) that there is a walk of length $\ell$ from a vertex $u_0$ to a vertex $u_{\ell}$ if and only if there is an $\ell$-tuple $(a_0,a_1,\dots,a_{\ell-1}) \in \{1,2,\dots,k\}^{\ell}$ such that $u_0=\sum_{i=0}^{\ell-1} a_i(-k)^i +(-k)^{\ell}u_{\ell}$. In fact, this gives a one-one correspondence between $\{1,2,\dots,k\}^{\ell}$ and the walks of length $\ell$ starting at $u_0$. By Lemma \ref{lemmaW}, it now follows that $A^m=J+(-1)^mI$.

To conclude that the Hoffman polynomial of $\G$ is $x^m-(-1)^m$ (and not a proper divisor), we claim that the diameter $D$ is at least $m$. To show this, observe that $A^i$ has row sums $k^i$, hence every vertex has at most $1+k+\cdots+k^{m-1}$ vertices at distance at most $m-1$. However, $(A^2)_{00}=1$, which implies that vertex $0$ has at most $k+\cdots+k^{m-1}$ vertices at distance $m-1$. Because this number is smaller than $n$, this indeed shows that $D \geq m$, and hence the Hoffman polynomial has degree at least $m$, and therefore must be $x^m-(-1)^m$. Thus, the eigenvalues of $\G$ are $k$ and the complex $m^{\text{th}}$-roots of $(-1)^m$ (including $-1$ itself). It also follows that $\G$ is strongly $\ell$-walk-regular for $\ell \equiv 0$ and $1 \text{ (mod }m)$.

Note that using circulants, Lam \cite{Lam} constructed several other $01$-matrices $A$ such that $A^{m}$ is a linear combination of $J$ and $I$, however as digraphs they have loops (that is, $A$ does have ones on the diagonal). Note also that the particular case $m=2$ is not new: J{\o}rgensen \cite{J} constructed these strongly regular digraphs, unaware of Lam's work.
\end{example}

By taking the product of above examples and the all-ones matrix or their line digraphs, thus adding an eigenvalue $0$, we get examples with the maximum number of real eigenvalues.

\begin{example}\label{ex:6} Let $A$ be the adjacency matrix of a digraph such that $A^{m}$ is a linear combination of $I$ and $J$ (so with $\mu_m=\lambda_m$), as in Examples \ref{ex:4}. Then both its line digraph and the digraph with adjacency matrix $A \otimes J$ are strongly $(m+1)$-walk-regular. Besides the eigenvalues of $A$, they have an extra eigenvalue $0$, so starting from the digraphs of Example \ref{ex:4}, we obtain examples of strongly walk-regular digraphs with three and four real eigenvalues, which is the maximum number according to Proposition \ref{prop:3or4}. Moreover, these digraphs are strongly $\ell$-walk-regular for $\ell \equiv 1 \text{ (mod }m)$.
\end{example}

Examples \ref{ex:4} and \ref{ex:6} show the typical behavior of the digraphs with a nonreal eigenvalue that are strongly $\ell$-walk-regular for infinitely many $\ell$, as we shall see in the next result.

\begin{theorem}\label{thm:nonreal} Let $\G$ be a strongly connected regular digraph with at least one nonreal eigenvalue. If $\G$ is strongly $\ell$-walk-regular for infinitely many $\ell$, then one of the following cases holds:

\begin{itemize}

\item[\rm(i)]  $\G$ is a doubly regular tournament;

\item[\rm(ii)] $\mul=\lal$ or $\mul=\nul$ for every $\ell$ for which $\G$ is strongly $\ell$-walk-regular;

\item[\rm(iii)] $\G$ is $k$-regular with four distinct eigenvalues
    $k,\rho e^{\pm i\varphi}$, and $\theta$, where $\theta$ is real and $0<|\theta|<\rho \leq k$.

\end{itemize}
\end{theorem}

\pf Let $\G$ be a $k$-regular strongly $\ell$-walk-regular digraph, with a nonreal eigenvalue $\rho e^{i \varphi}$, where $\rho>0$ and $\sin \varphi \neq 0$. Also its complex conjugate $\rho e^{-i \varphi}$ is an eigenvalue of $\G$, and by applying \eqref{eq:lswr1}, it follows that
\begin{equation}\label{eq:lswrnonreal}
\mul-\lal=-\frac{\sin \ell \varphi}{\sin \varphi}\rho^{\ell-1} \textrm{~~and~~} \mul-\nul=\frac{\sin (\ell-1) \varphi}{\sin \varphi}\rho^{\ell}.
\end{equation}

If $\G$ has only three distinct eigenvalues, then it is diagonalizable by Proposition \ref{multiplicities}, and hence it is strongly regular, and hence case (i) holds.

So we may assume from now on that $\G$ has at least four distinct eigenvalues, in particular let $\theta$ be another eigenvalue. By applying Corollary \ref{cor:three}, we find that
\begin{equation}\label{eq:theta}
\left(\frac{\theta}{\rho}\right)^{\ell} \sin \varphi =\frac{\theta}{\rho} \sin \ell \varphi - \sin (\ell-1)\varphi.
\end{equation}
Because this equation holds for infinitely many $\ell$ and its right hand side is bounded in $\ell$, it follows that $|\theta| \leq \rho$.

We now first consider the case that $\theta$ is nonreal, and aim to show that $\mul=\lal$ or $\mul=\nul$. By interchanging the roles of the eigenvalues in the above argument, it follows that $|\theta|=\rho$, and hence $\theta=\rho e^{i\varphi^*}$ for some $\varphi^*$ with $\sin \varphi^* \neq 0$ and $\cos \varphi^* \neq \cos \varphi$. Now \eqref{eq:theta} reduces to the following two equations:
\begin{align}\label{eq:1}
\sin \ell\varphi^* \sin \varphi&=\sin \varphi^* \sin \ell\varphi,\\
\cos \ell\varphi^* \sin \varphi&=\cos \varphi^* \sin \ell\varphi - \sin (\ell-1)\varphi.\label{eq:2}
\end{align}
Now let $r=\sin \ell \varphi /\sin \varphi$. If $r=0$, then by \eqref{eq:lswrnonreal}, we indeed have that $\mul=\lal$.
So we may assume that $r \neq 0$. Then \eqref{eq:1} and \eqref{eq:2} imply that
\begin{equation}\label{eq:r}r=\frac{\sin \ell \varphi^*}{\sin \varphi^*} = \frac{\cos \ell\varphi-\cos \ell\varphi^*}{\cos \varphi-\cos \varphi^*}.
\end{equation}

If  $\sin \varphi^*= \pm \sin \varphi$, then $\varphi^*= \pm \varphi +\pi$, and then \eqref{eq:r} implies that $\tan \ell \varphi=\tan \varphi$, so $(\ell-1) \varphi$ is a multiple of $\pi$. Now \eqref{eq:lswrnonreal} indeed implies that $\mul=\nul$.

If $\sin \varphi^* \neq \pm \sin \varphi$, then \eqref{eq:r} also implies that
\begin{equation*}\label{eq:r2}\frac{\cos \ell\varphi-\cos \ell\varphi^*}{\cos \varphi-\cos \varphi^*}=\frac{\sin \ell\varphi-\sin \ell\varphi^*}{\sin \varphi-\sin \varphi^*} =\frac{\sin \ell\varphi+\sin \ell\varphi^*}{\sin \varphi+\sin \varphi^*}.
\end{equation*}
Using these equations and sum-to-product trigonometric formulas, it follows that
\begin{equation*}\label{eq:tan}\tan\ell(\frac{\varphi+\varphi^*}{2})=\tan \frac{\varphi+\varphi^*}{2} \textrm{  and  }\tan\ell(\frac{\varphi-\varphi^*}{2})=\tan \frac{\varphi-\varphi^*}{2},
\end{equation*}
which again shows that $(\ell-1) \varphi$ is a multiple of $\pi$, and hence that $\mul=\nul$.

Next, assume that $\theta$ is real, with $|\theta|=\rho$. For the sake of readability, we will only consider the case that $\theta=\rho$. The case that $\theta=-\rho$ is similar, but a bit more technical. So we let $\theta=\rho$. Now \eqref{eq:theta} reduces to
\begin{equation*}\label{eq:thetarho}
\sin\varphi=\sin\ell\varphi-\sin(\ell-1)\varphi,
\end{equation*}
which is equivalent to
\begin{equation*}\label{eq:thetarho}
\sin\frac{\varphi}{2}\cos\frac{\varphi}{2}=\sin\frac{\varphi}{2}\cos\frac{(2\ell-1)\varphi}{2}.
\end{equation*}
Because $\sin\frac{\varphi}{2}\neq 0$, this implies that $\ell\varphi$ or $(\ell-1)\varphi$ is a multiple of $2\pi$, and hence by \eqref{eq:lswrnonreal}, we indeed have that $\mul=\lal$ or $\mul=\nul$.

Also the case that $\theta=0$ gives case (ii). Indeed, if $\theta=0$, then \eqref{eq:theta} implies that $\sin (\ell-1) \varphi=0$, and hence that $\mul=\nul$ by \eqref{eq:lswrnonreal}.

Finally, assume that we are not in any of the above cases and hence that $\theta$ is real, with $0<|\theta|<\rho$. We can rewrite \eqref{eq:theta} as
\begin{equation}\label{eq:cot}
(\tfrac{\theta}{\rho})^{\ell} =(\tfrac{\theta}{\rho}-\cos \varphi) \sin \ell \varphi / \sin \varphi+ \cos \ell\varphi.
\end{equation}
This implies that if $\{\ell_i\}_{i=1}^{\infty}$ is an increasing sequence such that $\G$ is strongly $\ell_i$-walk-regular for every $i$,
and $\lim_{i\rightarrow\infty}\sin \ell_i\varphi=0$, then also $\lim_{i\rightarrow\infty}\cos \ell_i\varphi=0$, which is a contradiction.
Thus, there is an increasing sequence $\{\ell_i\}_{i=1}^{\infty}$ such that $\G$ is strongly $\ell_i$-walk-regular and for which $|\sin \ell_i \varphi|>\delta$ for every $i$ and some $\delta>0$. Now it follows from \eqref{eq:cot} that $\lim_{i\rightarrow\infty}\cot \ell_i\varphi=(\cos \varphi-\theta/\rho)/\sin \varphi$. Thus, $\theta$ is determined by this equation, and hence we have the final case of the statement. \epf
\begin{corollary}\label{cor:inf} Let $\ell^*>1$ and let $\G$ be a strongly connected regular digraph that is strongly $\ell^*$-walk-regular with $\mu_{\ell^*}=\lambda_{\ell^*}$ or $\mu_{\ell^*}=\nu_{\ell^*}$, such that $\G$ is not strongly regular. Then $\G$ is strongly $\ell$-walk-regular for infinitely many $\ell$, and $\mul=\lal$ or $\mul=\nul$ for each $\ell$ such that $\G$ is strongly $\ell$-walk-regular. Moreover, if $m$ is the smallest integer such that $\G$ is strongly $m$-walk-regular with $m>1$, then the following holds:
\begin{itemize}

\item[\rm(i)] If $\mu_m=\lambda_m$, then $\G$ is strongly $\ell$-walk-regular for every $\ell \equiv 0$ and $1 \textup{ (mod }m)$;

\item[\rm(ii)] If $\mu_m=\nu_m$, then $\G$ is strongly $\ell$-walk-regular for every $\ell \equiv 1 \textup{ (mod }m-1)$;

\end{itemize}
\end{corollary}

\pf If $\mu_{\ell^*}=\lambda_{\ell^*}$, then $A^{\ell^*} \in \langle I,J \rangle$, which clearly implies that $\G$ is strongly $\ell$-walk-regular for every $\ell \equiv 0$ and $1 \textup{ (mod }\ell^*)$. If $\mu_{\ell^*}=\nu_{\ell^*}$, then $A^{\ell^*} \in \langle A,J \rangle$, which implies that $\G$ is strongly $\ell$-walk-regular for every $\ell \equiv 1 \textup{ (mod }\ell^*-1)$. Moreover, in both cases, each nonzero eigenvalue of $\G$ must have the same absolute value and has multiplicity one in the Hoffman polynomial, so $\G$ is diagonalizable. Hence, if $\G$ has nonreal eigenvalues, then by Theorem \ref{thm:nonreal} the result follows. If $\G$ has real eigenvalues only, then it must have three distinct eigenvalues besides the degree, and these must be $0$ and $\pm \rho$ for some $\rho$. In this case, it follows that $m=3$ and $\G$ is strongly $\ell$-walk-regular for every odd $\ell$, with $\mul=\nul$, see also Section \ref{sec:realdiagonalizable}.
\epf

Because the sum of the two nontrivial eigenvalues of a doubly regular tournament is $-1$, we can apply the same construction methods as in Examples \ref{ex:2} and \ref{ex:3}. In this case, we obtain examples with eigenvalues as in the final case of Theorem \ref{thm:nonreal}.

\begin{example}\label{ex:tournamentproduct} Consider a doubly regular tournament on $n$ vertices and let $A$ be its adjacency matrix. As mentioned before, its eigenvalues are $k=\frac12(n-1)$ and $-\frac12 \pm \frac12 \sqrt{-n}$. Then the digraph with adjacency matrix $J_{3n}-I-A \otimes J_3$ has distinct eigenvalues $\frac{3n+1}2, \frac12 \pm \frac32 \sqrt{-n}$, and $-1$. It thus follows from the Hoffman polynomial that it is strongly $3$-walk-regular.
\end{example}

We suspect however that the digraphs of this example are only strongly
$\ell$-walk-regular for $\ell=3$. It would therefore be interesting to find examples for case (iii) of Theorem \ref{thm:nonreal} that are strongly $\ell$-walk-regular for infinitely many $\ell$, or to show that no such examples exist.

\end{document}